\documentclass{article}

\usepackage{arxiv}
\usepackage[round]{natbib} 
\setcitestyle{authoryear,open={(},close={)}} 

\usepackage[utf8]{inputenc} 
\usepackage[T1]{fontenc}    
\usepackage{hyperref}       
\usepackage{url}            
\usepackage{booktabs}       
\usepackage{amsfonts}       
\usepackage{nicefrac}       
\usepackage{microtype}      
\usepackage{lipsum}
\usepackage{graphicx}
\usepackage{amsmath}
\usepackage{algorithm}
\usepackage{algpseudocode}
\usepackage{subcaption}
\usepackage{caption}
\graphicspath{ {./images/} }

\title{Automatic Calibration of Mesoscopic Traffic Simulation Using Vehicle Trajectory Data}

\author{
 Ran Sun \\
Department of Civil and Environmental Engineering, University of Michigan\\
Ann Arbor, MI 48109\\
  ransun@umich.edu\\
   \And
Zihao Wang \\
    Department of Civil and Environmental Engineering, University of Michigan\\
Ann Arbor, MI 48109\\
  zihaooo@umich.edu \\
  \And
Xingmin Wang \\
    Department of Civil and Environmental Engineering, University of Michigan\\
Ann Arbor, MI 48109\\
  xingminw@umich.edu \\
  \And
 Henry X. Liu*\\ 
  Department of Civil and Environmental Engineering, 
  University of Michigan\\
Ann Arbor, MI 48109\\
  henryliu@umich.edu
}

\begin{document}
\maketitle
\begin{abstract}
Traffic simulation models have long been popular in modern traffic planning and operation applications. Efficient calibration of simulation models is usually a crucial step in a simulation study. However, traditional calibration procedures are often resource-intensive and time-consuming, limiting the broader adoption of simulation models.  In this study, a vehicle trajectory-based automatic calibration framework for mesoscopic traffic simulation is proposed. The framework incorporates behavior models from both the demand and the supply sides of a traffic network. An optimization-based network flow estimation model is designed for demand and route choice calibration. Dimensionality reduction techniques are incorporated to define the zoning system and the path choice set. A stochastic approximation model is established for capacity and driving behavior parameter calibration. The applicability and performance of the calibration framework are demonstrated through a case study for the City of Birmingham network in Michigan.
\end{abstract}

\keywords{ Traffic simulation calibration\and Vehicle trajectory data \and  OD estimation\and Simulation-based optimization}

\newpage

\section{Introduction}
Traffic simulation plays a crucial role in many transportation operations and planning applications. Large-scale, high-fidelity simulation modeling has enabled cost-effective and flexible assessment of traffic design and management options. These simulations are particularly valuable due to their comprehensive outcomes, such as evaluating the impacts of traffic signal control, pricing and toll strategies, safety and risk prevention designs, and the introduction of intelligent transportation systems (ITS). Despite their wide applications, timely and adequate calibration of these simulation models for daily operations can be challenging due to a lack of data, improper use of calibration tools, or constraints in resources and time.

Significant research has focused on the calibration of traffic simulation models. On the supply side, parameters such as link and node capacity and driving behavior influence traffic flow within the network. These parameters are adjusted manually or automatically to ensure model outputs align with observed data. 
To calibrate driving behavior and capacity parameters, it is important to include data from a wide range of traffic conditions (i.e. free flow and congested situations) and from a wide network coverage. Due to the limitation in the real world, many studies only focus on a subset of parameters to adjust depending on the specific traffic simulation model to calibration. 
Various methods have been used to calibrate the traffic simulator, including sensitive analysis approaches in
\cite{ciuffo2014sensitivity}, step-wise multistage calibration procedure in \cite{chu2003calibration}, stochastic approximation method in \cite{paz2015calibration}. The performance is often evaluated based on various goodness-of-fit measures for traffic volume, space mean speed, density, and queue length.

On the demand side, Origin-Destination (OD) demand and route choice serve as crucial inputs to most large-scale traffic simulation models. The route choice component is sometimes jointly considered with the calibration of OD flow. Because the route choice bridges the OD flow, path flow, and link flow in the simulation network. The goal of OD demand calibration or estimation is to estimate a trip table output, often in terms of trip origin, trip destination, departure time, and sometimes the chosen routes.It is generally considered to be an under-determined problem because the number of OD pairs is significantly larger than the number of observable links or paths, leading to multiple possible OD demand patterns that can produce the same system state.  Many studies approach the problem from statistical inference of the OD flow, including Bayesian estimation \citep{hazelton2008statistical} and statistical properties of the OD flow \citep{shao2014estimation}. Scalability for estimation in large-scale networks has also been addressed using stochastic programming-based OD estimation \citep{sun2024stochastic} and metamodel approximation methods \citep{osorio2019high}. Even with the opportunities shown in the literature, most of the methods still rely on link counts to estimate the system demand, which fundamentally suffers from the ill-possess of the problem structure. 

Even with the ongoing research activities, 
calibrating large-scale traffic simulation models is a still resource-intensive and time-consuming process that often requires substantial manual intervention and adjustments from data collection to trial-and-error verification and validation \citep{dowling2004traffic}.
Extensive manual work would be needed to fine-tune various model parameters to ensure that the simulation accurately reflects real-world traffic conditions. The scale and complexity of transportation systems operations can result in significant costs and extended durations for the calibration process. Additionally, reliance on manual correction increases the potential for human error, further complicating and prolonging the calibration efforts. Consequently, developing more efficient and automated calibration techniques remains a critical area of research in traffic modeling.

Trajectory data opens the door for the development of an automatic calibration process. There are two unique advantages for trajectory data based calibration. First,  trajectory data directly provides partial path flow information.  Different from data collected from loop detectors, cameras or license place recognition, partial path flow from trajectory data provide wide coverage that could span the entire road network. Path-level information from trajectories such as origin and destination, route choice, and trip travel time can directly provide values for simulator parameter calibration. This fundamentally improves the underdetermined network flow estimation problem. Even with low penetration rate in space and time, path-level information from trajectory data would be beneficial for the understanding of the system operation. 
Second, trajectory based approach can systematically integrate the link based components. Among most existing studies, aggregated link-based detector data is still often used as the only data source for calibration. Trajectory based approach is able to combine the data from previously available detectors. Additionally, our previous work on link based penetration rate estimation at urban signalized intersections can be integrated \citep{wang2024traffic}.

In this paper, we establish an automatic calibration framework for mesoscopic traffic simulation using vehicle trajectory data. This framework systematically calibrates mesoscopic traffic simulation models by combining network coverage from path-level information and consistency regularization from link-level information from signalized intersection estimation models. Both demand and supply components are jointly calibrated to reflect real-world observations. Trajectory data is used as an approximation to the equilibrium state of network
for the 'replay' simulation of statistical average of the real-world. Dimensionality reduction and denoising techniques for demand clusters and path choice sets are incorporated to improve estimation quality and calibration accuracy. Unlike classical path generation methods restricted to a certain number of paths, our data-driven approach maximizes information retention from trajectory data.  The proposed trajectory-based calibration framework is demonstrated based on a case study of the City of Birmingham network. 
The automatic calibration framework would be able to extend to the calibration of any network with low-penetration vehicle trajectory data. 

The remainder of this paper is organized as follows. In the second section, a trajectory-based automatic mesoscopic simulation calibration framework is presented. The third section describes the detailed applications of calibrating the simulation for a case study in the City of Birmingham network. The final section presents the discussions and conclusions. 


\section{An Automatic Trajectory-based Mesoscopic Simulation Calibration Framework}

\begin{figure}[!htbp]
  \centering
  \includegraphics[width=0.8\textwidth]{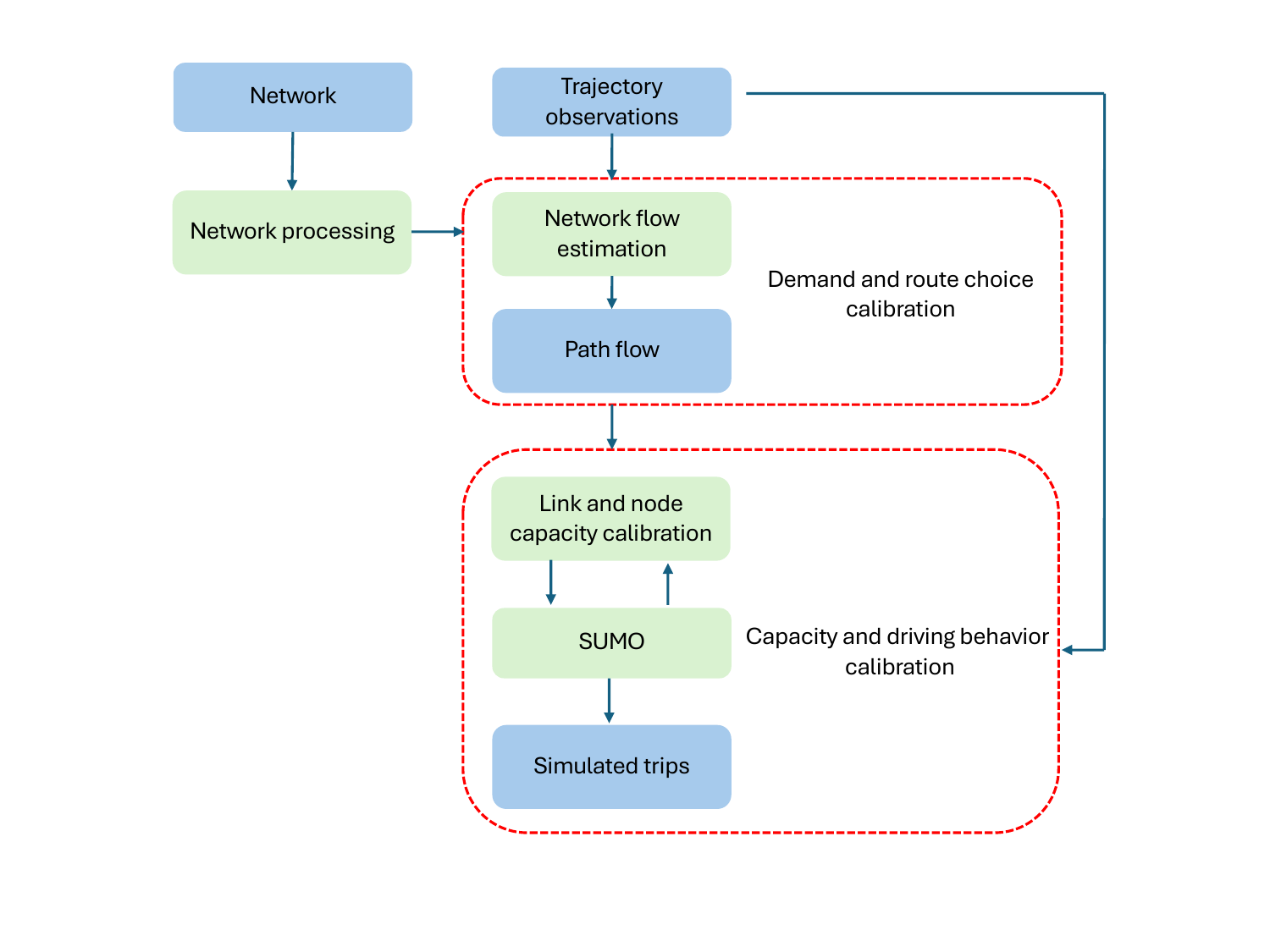}
  \caption{An Automatic Calibration framework}\label{fig:framework}
\end{figure}

In this section, we present our proposed framework for automatic trajectory-based simulation calibration. The framework jointly calibrates the demand and route choice, and capacity and driving behavior parameters for the network. In the first step, we process the network and trajectory data into a ready-to-use structure following dimensionality reduction and clustering techniques. In the second step, we formulate the network flow estimation as an estimation from partial observations problem. The goal is to reconstruct full-scale OD flow, path flow, and link flow for the entire network based on the given partial path flow data derived from trajectory data over the simulation time intervals. In the third step, based on the estimated network flow, a simulation based optimization is designed to optimize the link and node capacity and driving behavior parameters for the simulation model. The goal is to match the performance metrics using simulated trips and observed trajectory samples in the network. The key steps of the framework is presented in Figure \ref{fig:framework}.

\subsection{Demand and route choice calibration: a network flow estimation model}
OD demand, also referred to as OD flow in this study, determines the number of trips between each origin and destination zone pair during a specific departure time interval. 
In order to prepare a comprehensive set of simulation inputs, one needs to start with estimating the trip table. 
A crucial consideration is that the calibrated simulation must reflect the system throughput in terms of the number of trips the network can support across various time intervals. Based on the methodology in \citep{wang2024traffic}, we utilize a node-level model to estimate the penetration rate at signalized intersections. The estimated penetration rate is to inversely compute the total trips, serving as the surrogate ground truth knowledge of the full-scale trips in the network across different time-of-day (TOD) intervals. Thus, a desired property of the network flow estimation model is that the total estimated flow should be consistent with the trajectory penetration rate.

The network flow estimation is formulated as an estimation from partial observations problem. 
The partial path flow comes from two sources: 1) Observations are available for only a subset of all paths in the network within each time interval, meaning many paths used in the real world are not covered by the trajectory data during that time; and 2) Even for paths with trajectory coverage, we lack the ground truth for the full-scale path flow, i.e. the trajectory sampling rate for each specific path. Consider a network with $N$ OD pairs, $M$ paths and $E$ links. Denote OD flow as $\boldsymbol{x}\in \mathbb{R}^N$, path flow as $\boldsymbol{y}\in\mathbb{R}^M$ and link flow as $\boldsymbol{z}\in\mathbb{R}^E$. We now present the proposed network flow estimation model as in optimization program \ref{eq:flow_opt} below,

\begin{flalign}
\min_{\boldsymbol{x}, \boldsymbol{y}, \boldsymbol{z}} & \quad \| \boldsymbol{1}^\top \boldsymbol{x} - \tilde{x}\|_2^2 + \gamma \|\boldsymbol{y} - \boldsymbol{G} \boldsymbol{x}\|_2^2 + \rho \| \boldsymbol{W} (\boldsymbol{z} - \boldsymbol{\tilde{z}}) \|_2^2 \\
\text{s.t.} & \quad \boldsymbol{x} = \boldsymbol{\Phi} \boldsymbol{y} \quad \text{(OD-Path incidence)} \\
& \quad \boldsymbol{z} = \boldsymbol{\Omega} \boldsymbol{y} \quad \text{(Link-Path incidence)} \\
& \quad x_i \in [0, \alpha_i], \quad i = 1, 2, \ldots, N \quad \text{(OD box constraint)} \\
& \quad y_i \in [0, \beta_i], \quad i = 1, 2, \ldots, M \quad \text{(Path box constraint)} \\
& \quad z_i \in [0, C_i], \quad i = 1, 2, \ldots, E \quad \text{(Link box constraint, capacity)}.
\end{flalign}
\label{eq:flow_opt}

The network flow estimation is formulated as a quadratic program (QP). The problem is solved independently for each TOD time interval. The goal is to estimate OD flow, path flow and link flow based on total full-scale trip volume, network topology, network operation knowledge, and partial link flow estimations obtained for a subset of links connecting signalized intersections. The objective function aims to minimize the sum of three error components: the deviation of the total estimated trips in the region from its prior full-scale OD flow, the correspondence of the estimated path flows and OD flow following an approximated equilibrium state, and the deviation of the estimated link flows from field link flowestimates obtained from signalized intersections estimation model. In the first term, the summation of OD flow $\boldsymbol{1^\top x}$ should be close to total trips $\tilde{x}\in\mathbb{R}$, with $\boldsymbol{1}$ being an all-one summation vector. In the second term, the path flow $\boldsymbol{y}$ should align with the OD-path assignment $\boldsymbol{Gx}$. The $\boldsymbol{G}\in\mathbb{R}^{M\times N}$ is an exogenous traffic assignment matrix translating OD flow to path flow. In this study, we use $\boldsymbol{G}$ as an approximation to the equilibrium state on an average case for each TOD interval. The last term represents the weighted link flow error. The diagonal weight matrix $\boldsymbol{W}\in\mathbb{R}^{E\times E}$ encodes the indices of the available link flow locations in the network and their confidence rate. The estimated link flow $\boldsymbol{z}\in \mathbb{R}^E$ should match the zero-padded link flow prior  $\boldsymbol{\tilde{z}}\in \mathbb{R}^E$. The link flow prior is estimated based on trajectory data following the previous work in \citep{wang2024traffic}. The weighting parameters $\gamma$ and $\rho$ control the balance for the error terms. 

Two types of information are included in the constraints. First, OD-path incidence $\boldsymbol{\Phi}\in\mathbb{R}^{N\times M}$ represents the correspondence between each OD pair and its path set. The link-path incidence $\boldsymbol{\Omega}\in\mathbb{R}^{E\times M}$ represents the set of links belonging to each path. Second, the box constraints for OD flow, path flow, and link flow control the non-negativity and ensure the realistic magnitude of these components, respectively.

Given the network topology and its corresponding link-path and OD-path incidence matrices, we can reformulate the problem into an equivalent QP with only path flow $\boldsymbol{y}$ as the decision variable, as shown in program \ref{eq:optimization_problem_path}.

\begin{flalign}
\min_{\boldsymbol{y}} & \quad \|\boldsymbol{1}^\top \boldsymbol{\Phi} \boldsymbol{y} - \boldsymbol{\tilde{x}}\|_2^2 + \gamma \|\boldsymbol{y} - \boldsymbol{G} \boldsymbol{\Phi} \boldsymbol{y}\|_2^2 + \rho \| \boldsymbol{W} (\boldsymbol{\Omega} \boldsymbol{y} - \boldsymbol{\tilde{z}}) \|_2^2 & \\
\text{s.t.} & \quad (\boldsymbol{\Phi} \boldsymbol{y})_i \in [0, \alpha_i], \quad i = 1, 2, \ldots, N & \\
& \quad y_i \in [0, \beta_i], \quad i = 1, 2, \ldots, M & \\
& \quad (\boldsymbol{\Omega} \boldsymbol{y})_i \in [0, C_i], \quad i = 1, 2, \ldots, E. &
\end{flalign}
\label{eq:optimization_problem_path}

The problem is reduced to determining the optimal path flow that matches the total OD scale and estimated link flow, with constraints ensuring the path flows remain within specified bounds. Popular ways to solve QP problems include interior-point methods, active-set methods, and gradient-based methods. Interior-point methods efficiently handle large-scale QPs by traversing the interior of the feasible region. Active-set methods focus on the constraints that are active at the optimal solution, iteratively adjusting the active set to find the optimum. Gradient-based methods, including conjugate gradient and quasi-Newton methods, leverage the problem's gradient and Hessian information to iteratively improve the solution. In this study, we use an ADMM-based first-order method implemented in OSQP due to its preferred robustness and scalability \citep{osqp}.

Note that in the network flow estimation, the problem is treated it as a continuous optimization program.
However, in SUMO traffic simulation, a discrete trip table input and its corresponding route choice are preferred, because we could have better control for trip level information. Therefore, based on the continuous path flow values directly obtained from the optimization program, we post-process the results by rounding the path flows to the nearest integer to ensure feasibility for traffic simulation.

\subsection{Capacity and driving behavior calibration: a stochastic approximation model}
In this section, we present the second main component in the estimation framework, capacity and driving parameter calibration. Specifically, we focus o SUMO parameters including 1) mesoscopic traffic model parameters for the longitudinal model and junction model; 2) departure time-related trip loading parameters; 3) vehicle attributes; and 4) re-routing parameters. Collectively, these parameters define the link and node capacity for the simulation network. The calibration problem can be summarized as minimizing the difference in simulated trip travel time and observed trip travel time, as following,

\begin{flalign}
\min_{\boldsymbol{\theta}} & \quad     L(\boldsymbol{\theta}) = \frac{1}{K}\sum_{k=1}^K (S_k(\boldsymbol{\theta};\boldsymbol{y}) - O_k)^2,\\
\text{s.t.} & \quad \boldsymbol{\theta}\in \boldsymbol{\Theta}.
\end{flalign}
\label{eq:spsa}

$S(\cdot)$ represents the simulation-based function for travel time of simulated trip $k$ from SUMO, $O_k$ is observed trip travel time from trajectory for trip $k$. Path flow $\boldsymbol{y}$ is obtained from the flow estimation model \ref{eq:flow_opt}. The error function $L(\cdot)$ measures the discrepancy in travel time that we aim to minimize. The constraint $\boldsymbol{\Theta}$ defines the feasible values for the simulation parameters to be calibrated.  

The main challenge is that the simulation outcome in terms of trip statistics is a non-convex function w.r.t. the simulation parameters due to the involvement of the simulation-based function $S(\cdot)$. Traditional convex optimization algorithms cannot be directly applied. Moreover, the simulation-based optimization is computationally intensive, requiring efficient methods for a tractable global search.  This is because a mathematically expensive simulation is called every iteration when the parameters are updated. Simultaneous Perturbation Stochastic Approximation (SPSA) is an iterative optimization method specifically designed for optimizing complex problems where the objective function is noisy or expensive to evaluate \citep{spall1992multivariate}. SPSA replies on the first-order derivative estimates of the gradient of the objective function using random perturbations. It estimates the gradient with two evaluations of the objective function per iteration, regardless of the dimensionality of the decision variable. This characteristic makes SPSA specifically suitable for high-dimensional non-convex problems. The algorithm for our calibration problem is summarized in Algorithm \ref{algo:spsa}. In each iteration for SUMO simulation, we use the first TOD interval as the warm-up stage and the last TOD interval as the cool-down stage. The update and performance evaluation are conducted based on the main TOD intervals from 7 am to 10 pm.

\begin{algorithm}
\caption{SPSA for Mesoscopic Simulation Calibration}
\begin{algorithmic}[1]
\State \textbf{Initialization:}
\State Choose initial guess $\boldsymbol{\theta}_0$
\State Choose sequences $a_k$, $c_k$, and perturbation vector $\boldsymbol{\Delta}_k$
\State Set iteration counter $k = 0$

\While {Stopping criteria: changes in total travel time > $\epsilon$}
    \State Generate random perturbation vector $\boldsymbol{\Delta}_k$
    \State Ensure perturbed points are in the feasible region $\boldsymbol{\Theta}$
    \State Evaluate objective function at perturbed points:
    \[
    L(\boldsymbol{\theta}_k + c_k \boldsymbol{\Delta}_k) \quad \text{and} \quad L(\boldsymbol{\theta}_k - c_k \boldsymbol{\Delta}_k)
    \]
    \State Estimate gradient $g_k$ of objective function:
    \[
    g_k = \frac{L(\boldsymbol{\theta}_k + c_k \boldsymbol{\Delta}_k) - L(\boldsymbol{\theta}_k - c_k \boldsymbol{\Delta}_k)}{2c_k} \boldsymbol{\Delta}_k^{-1}
    \]
    \State Update parameters:
    \[
    \boldsymbol{\theta}_{k+1} = \boldsymbol{\theta}_k - a_k \boldsymbol{g}_k
    \]
    \State Increment iteration counter $k = k + 1$
\EndWhile
\State \textbf{Output:} Optimal parameters $\boldsymbol{\theta}_{opt} = \boldsymbol{\theta}_k$
\label{algo:spsa}
\end{algorithmic}
\end{algorithm}

\section{Case Study: Calibration of Birmingham Network}
The network for case study is a congested mid-scale network in the City of Birmingham, Michigan. We obtained the trajectory from an Original Equipment Manufacturer (OEM) with the penetration rate estimated as approximately $7.5\%$ based on a previous study \citep{wang2024traffic}. This translates to 0.28 million trips in total for an average weekday in the city. The link flow prior is based on the estimation of 12 links along the Adams Road corridor for TOD intervals 1 to 4. The outcomes for the automatic calibration framework include the simulated trips for one entire weekday, calibrated demand and supply parameters, and various performance measures.

\subsection{Network processing and trajectory map-matching }
\begin{figure}[!htbp]
  \centering
  \begin{subfigure}[b]{0.4\textwidth}
    \centering
    \includegraphics[width=\textwidth]{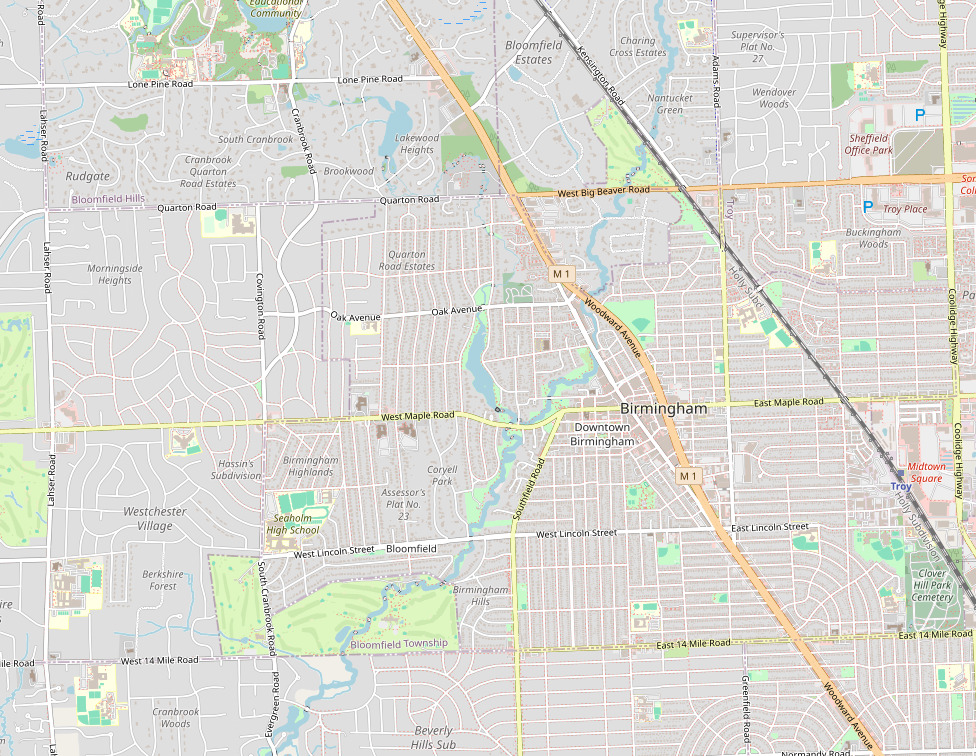}
    \caption{Geospatial data from OSM}\label{fig:network_subfig1}
  \end{subfigure}
  \hfill
  \begin{subfigure}[b]{0.36\textwidth}
    \centering
\includegraphics[width=1\textwidth]{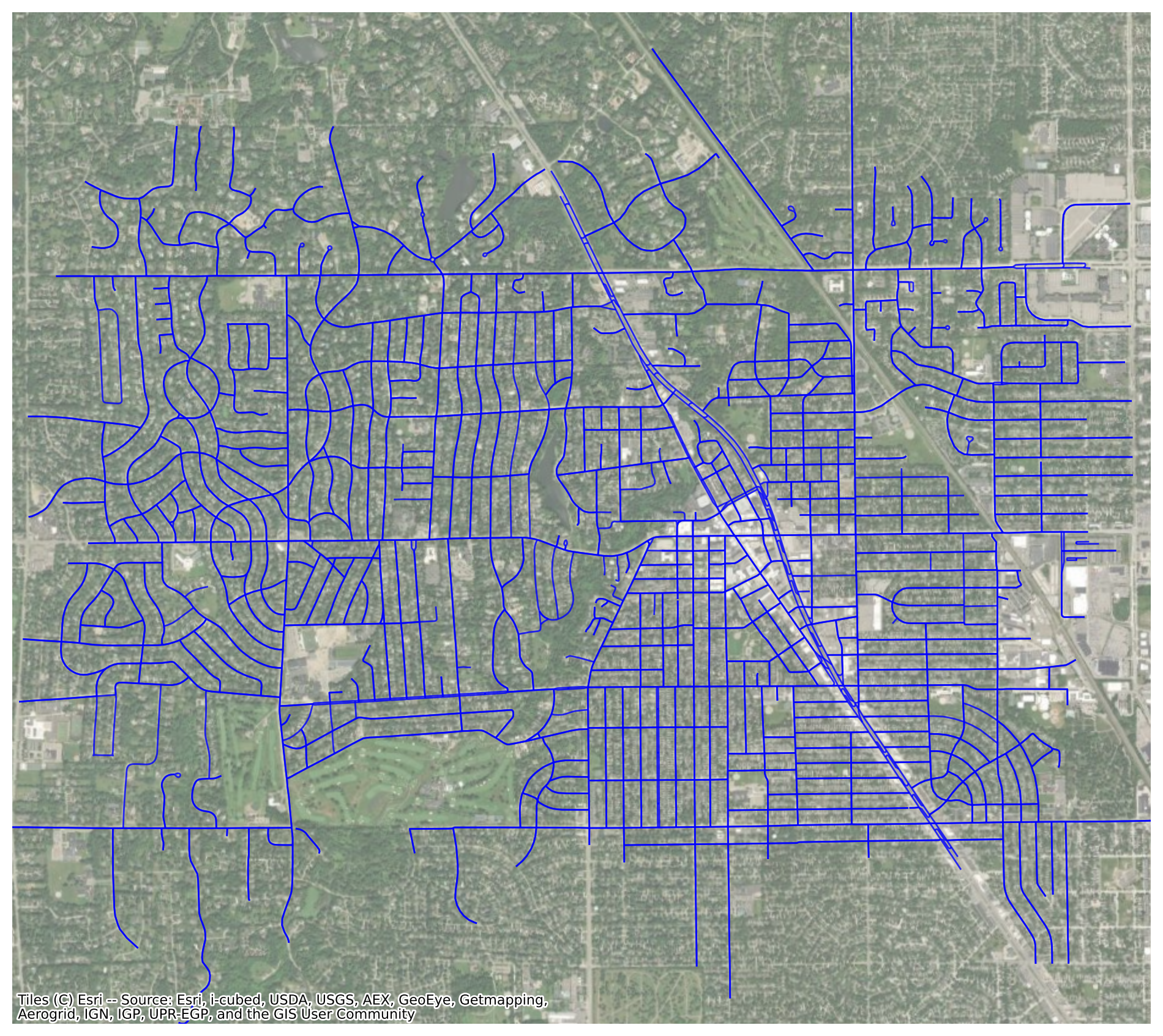}
    \caption{Road network}\label{fig:network_subfig2}
  \end{subfigure}
  \caption{Birmingham road network}\label{fig:network}
\end{figure}

The original geospatial data for the Birmingham network was obtained from OSM\footnote{https://www.openstreetmap.org/}, as shown in Figure \ref{fig:network_subfig1}. OSM is a widely used open-source mapping platform that provides rich geographical data, including road network geometry and associated speed limits.  The raw OSM network was then converted to the road network format ready for trajectory map-matching. The processed network includes 3066 links and 8443 nodes, as shown in Figure \ref{fig:network_subfig2}. The Signal Phase and Timing data (SPaT) provided by the local traffic management center was integrated into the network. 

The raw trajectory data consists of timestamped location points that provide a detailed record of vehicle movements over time. The data was pre-processed to remove any erroneous records that could result from GPS inaccuracies or data logging errors. The processing procedure is detailed in \citep{wang2023trajectory}.

Despite the convenient process of integrating OSM data, the speed limit data in OSM can often be outdated or incorrect due to infrequent updates or user input errors. With the availability of vehicle trajectory data after map matching, we are able to leverage vehicle trajectory speed observations to compute more accurate estimates of speed limits. Specifically, we compute the $80$ percentile of observed free-flow speed for each link that is covered by trajectory data. The updated speed limits are integrated back into the network and are subsequently used in traffic simulations and other analyses to ensure more reliable results.

In addition, there exist abnormal trips with GPS errors and extremely slow trips, such as idling off-street or constant GPS drifting. These trips mostly occur off the road network, for example, in parking spaces or residential personal spaces, and do not directly contribute to network congestion. Therefore, abnormal trips with average speeds less than 5 mph or greater than 100 mph are removed from the simulation calibration and performance evaluation process.

\subsection{OD demand clustering}
Defining the TAZ system is a fundamental prerequisite for the estimation process. In practice, TAZs systems are often developed by local Metropolitan Planning Organizations (MPOs) to facilitate planning and operation applications. However, the default TAZ system might not reflect the information on realized trips in the region. Trips may be unevenly distributed across different TAZs. TAZs may be coarse in certain regions where we need fine-resolution analysis, for example, signal control design in dense grid networks and search for parking in central business areas. More importantly, the TAZ system controls the dimension $N$ of the OD flow to be estimated for the network. Finding a flexible TAZ system would act as a dimensionality reduction processing and benefit the calibration steps. 

\begin{figure}[!htbp]
  \centering
  \includegraphics[width=0.6\textwidth]{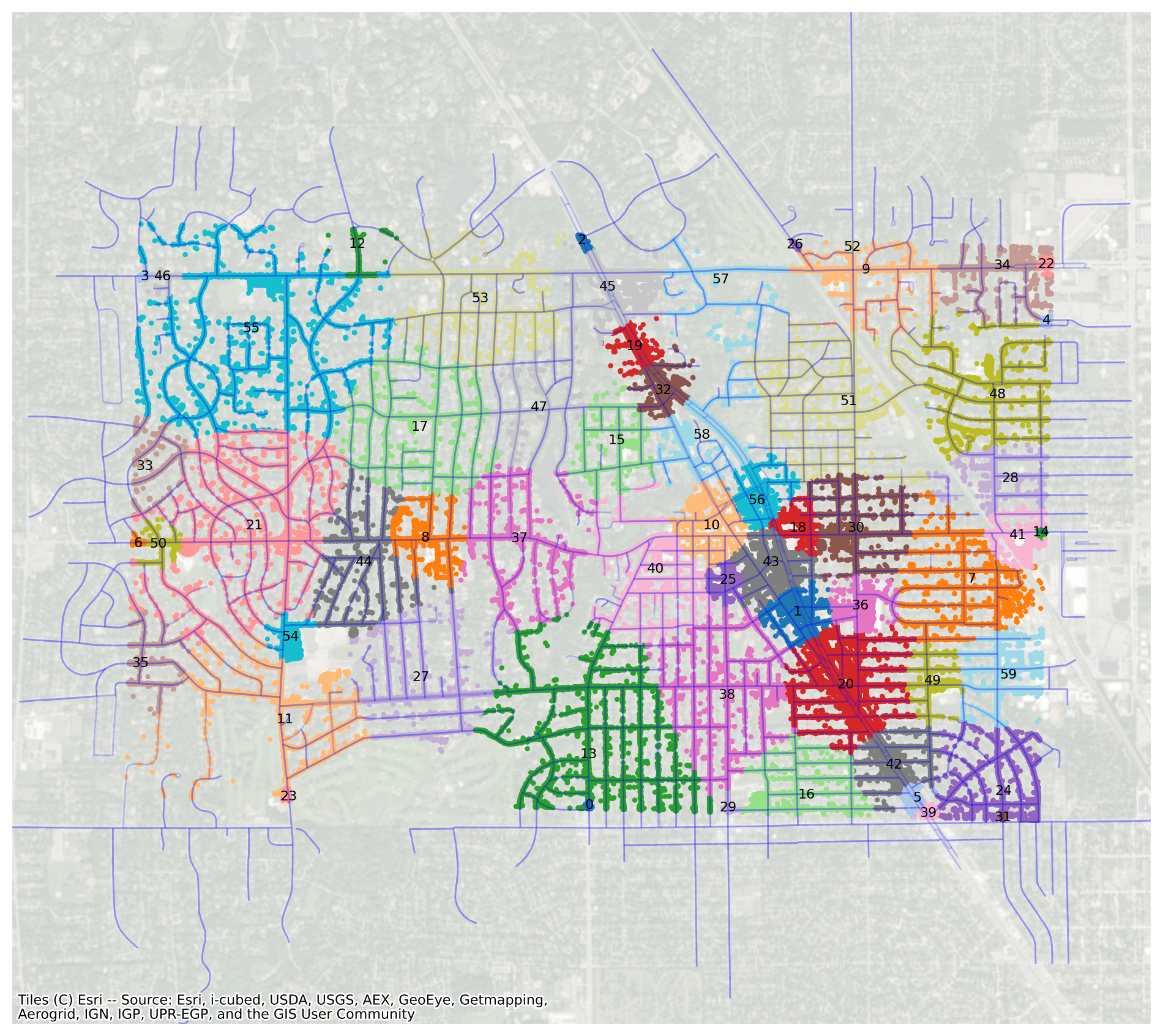}
  \caption{Demand clustering}\label{fig:od_clustering}
\end{figure}

Based on trajectory data, we design a data-driven approach that groups the origin and destination locations in the region. Specifically, the latitude and longitude for each origin point or destination point is treated as one data sample for the clustering. We want to cluster all these points based on the local density across the region. In this study, we utilize a Gaussian mixture model for demand clustering in the study region. 
The resulting zoning system consists of 60 TAZs, as shown in Figure \ref{fig:od_clustering}. The resulting point-based clusters are then translated to network link-based clusters from a 'majority-vote' processing for each link. 

\subsection{Path set generation from path clustering }

To implement the flow estimation model, one needs to derive the path set from observed vehicle trajectory data. This step is important for defining the path flow domain and establishing OD-Path-Link correspondence. However, there are generally three challenges in constructing the path set. First, the dimensionality of the path set is very high, even for a moderate-sized network. Enumerating all feasible paths in the network is extremely computationally expensive. In practice, often a simplified path set is used for network flow studies, for example K-shortest path used in \citep{van2005path} and the sub-networks designed in \citep{dial2006path}.  Second, the paths often have high overlap. Two paths connecting one OD pair may only deviate in a few links. This would cause the violation of the Independence from Irrelevant Alternatives (IIA) property in route choice modeling. Accordingly, several route choice models are proposed, such as the 
c-logit model \citep{cascetta1996modified} and path-size logit \citep{ben1999discrete}. Third, the directly observed trajectory can be noisy. Uncommon behavior can cause the failure to capture the utility and route choice. Examples include taxi pick-up, parking search, and detours due to misguided routes. 

Our path flow estimation model would start from a k-shortest path enumeration for the network. This means a path enumeration algorithm is needed to generate the candidate path choice set. 
In our study, with the availability of trajectory data, we are able to derive the realized path from real-world data. From trajectory observations during multiple days, it is reasonable to assume that the realized path can serve as a reduced candidate path set. This mitigates the effort of conducting complex network search. In addition, the path clustering technique would help us identify representative paths used for all the OD pairs in the network. This also serves as a denoising and dimensionality reduction step to help us keep the problem more trackable. However, it is worth noting that this assumption means the feasible path that is not covered by any trajectory over a month of data is not considered in the path set. In the future, we hope to improve this assumption by combining the path enumeration and trajectory based realized paths for better efficiency and complexity trade-off. 

Based on the network topology, we can obtain the link-path incidence matrix $\boldsymbol{\Psi}\in\mathbb{R}^{E\times E}$. To further account for the heterogeneity in different links, we can augment the incidence matrix to the link length incidence matrix $\boldsymbol{\Psi}^L$, i.e. $\boldsymbol{\Psi}^L_{i,j} = l_j$ with $l_j$ being the length of link $j$. Next denote $\boldsymbol{\Psi}^L_i$ as the $i$th path in the network, with the elements being the lengths of the links along this path. Now we can define the similarity of path $i$ and path $i'$ using Jaccard similarity based on overlapping link lengths, as 

    \begin{equation}
        J(i, i') = \frac{|\boldsymbol{\Psi}^L_i \cap \boldsymbol{\Psi}^L_{i'}|}{|\boldsymbol{\Psi}^L_i \cup \boldsymbol{\Psi}^L_{i'}|}
    \end{equation}

From the Jaccard similarity as the pairwise similarity matrix for all candidate noisy paths, we apply a hierarchical clustering that groups the paths for each OD pair. As a result, the reduced path flow domain and reduced path-based incidence matrix are obtained, which are then used in the estimation of the network flow. The resulting path choice set is flexible in terms of the number of paths for each given OD pair, and this would reduce the unnecessary and unrealistic path chosen by typical path enumeration techniques. In addition, the generated representative path set serves as a dimensionality reduction process and contributes to the computation efficiency of the network flow estimation.

\subsection{Approximate assignment mapping}
\begin{figure}[!htbp]
  \centering
  \includegraphics[width=1\textwidth]{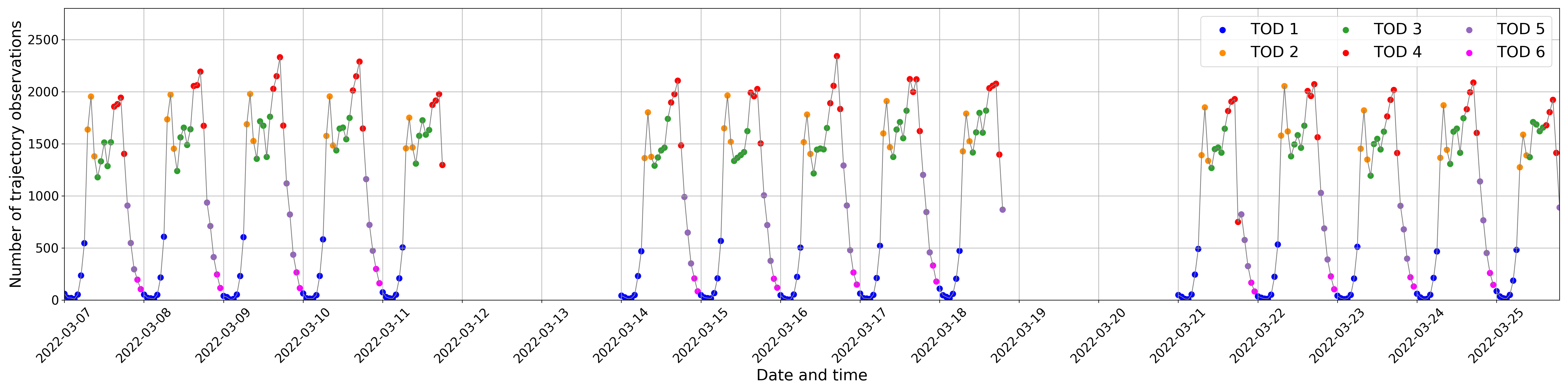}
  \caption{Trajectory observation counts}\label{fig:traj_time}
\end{figure}

\begin{figure}[!htbp]
  \centering
  \includegraphics[width=0.7\textwidth]{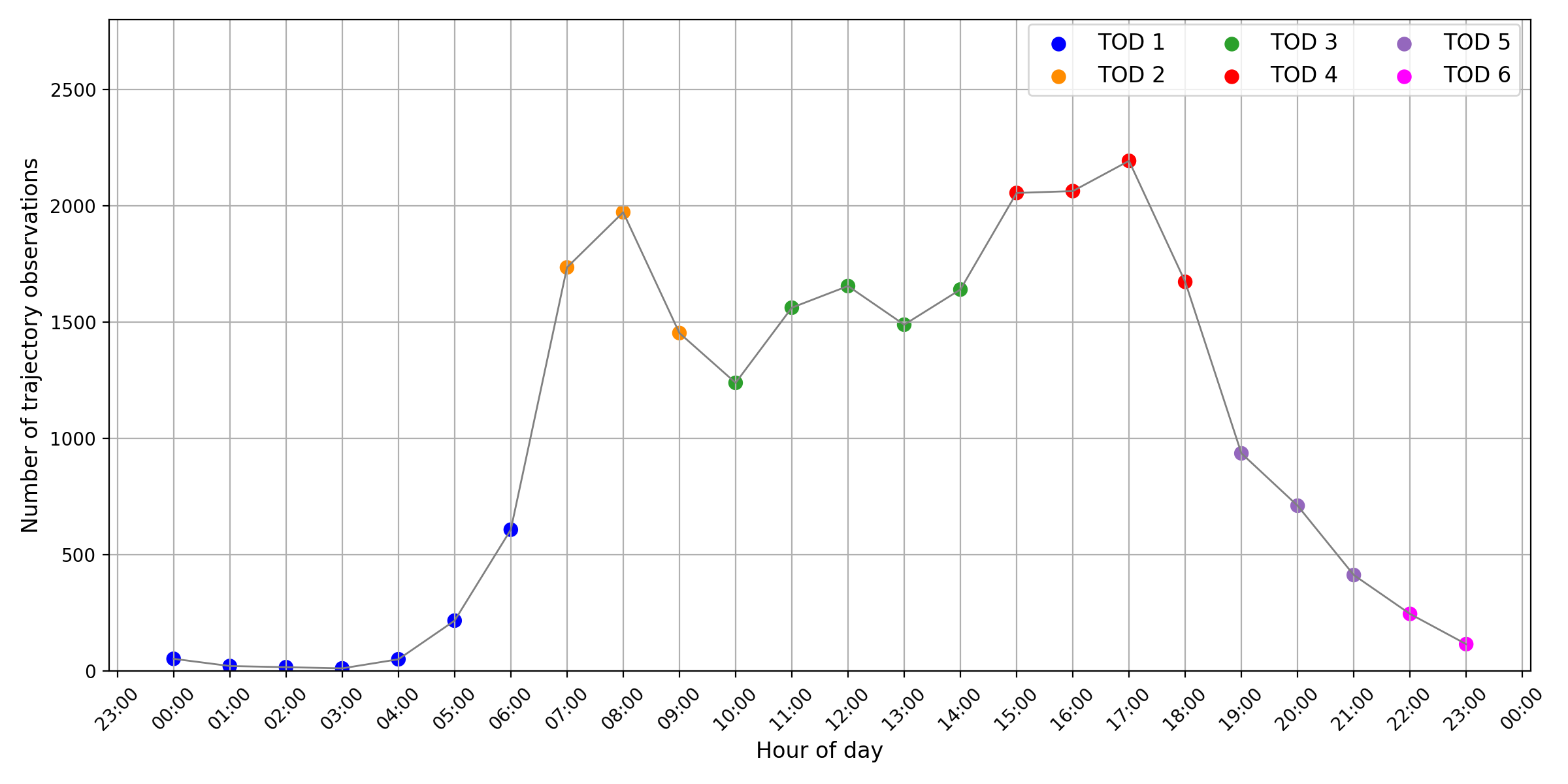}
  \caption{Trajectory observation counts over 2022-03-08}\label{fig:traj_time_day}
\end{figure}

In order to construct the traffic assignment mapping connecting OD flow and path flow. We start with a visualization of the total observed trajectory over time, shown in Figure \ref{fig:traj_time}. In the plot, each dot represents an hourly aggregation of all the trips in the network, with different color representing different TOD intervals. In Figure \ref{fig:traj_time_day}, we zoom into the total trips within one typical week day. 

Table \ref{tab:tod_intervals} presents the definition for TOD intervals. The TOD intervals are originally obtained from SPaT and it is found consistent with the temporal pattern of real-world trajectory. Within each TOD interval, total observed trips are relatively stationary across multiple days. This implies an assumption that the network state is stable within each TOD interval across multiple days. Therefore, we can use trajectory from multiple days for the same TOD interval to estimate the linear approximation to the equilibrium traffic assignment. As a result, six approximation matrices are computed for 6 TOD intervals. 

In this study, we assume the observed trajectory is a representative set of partial observations of the time-dependent equilibrium state of the road network. Therefore, we want to use SUMO simulation to best approximate the equilibrium state across all TOD intervals.

\begin{figure}[!htbp]
  \centering
  \includegraphics[width=1\textwidth]{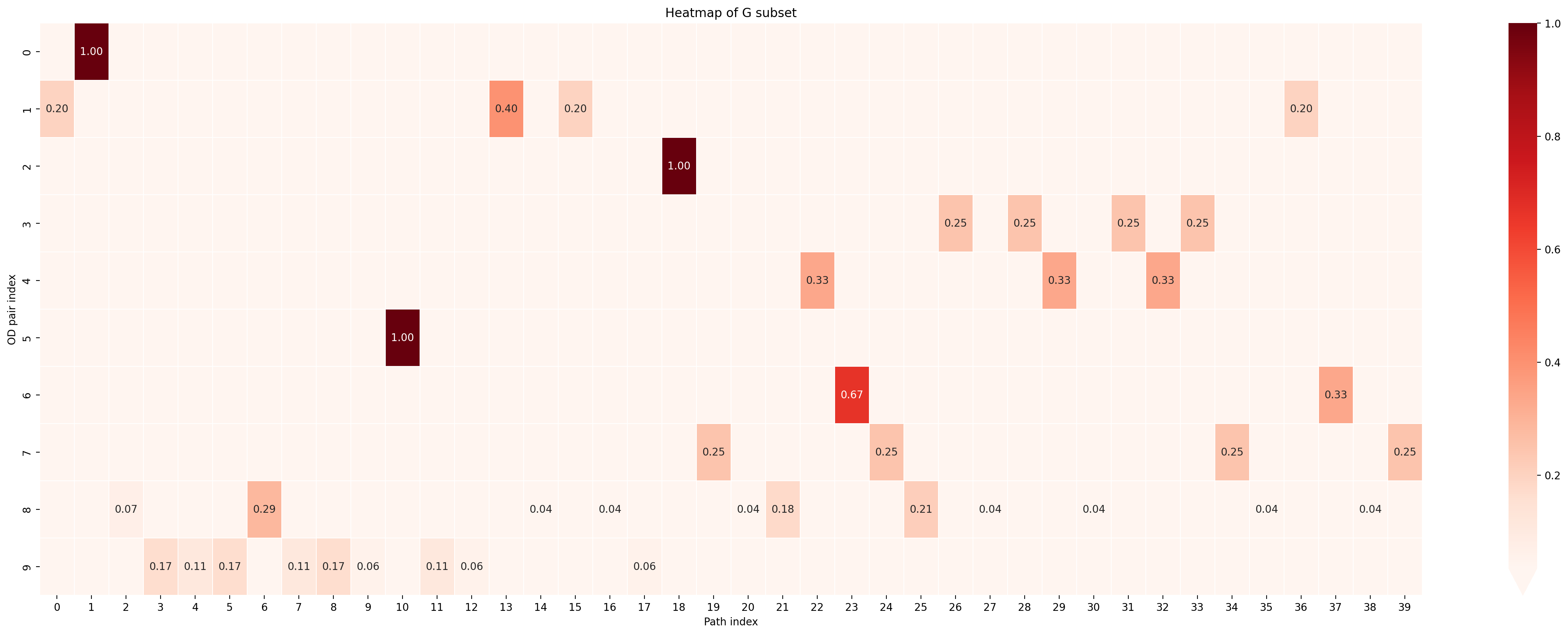}
  \caption{An example of approximate assignment mapping}\label{fig:G_example}
\end{figure}

Instead of calibrating a user equilibrium model, we use a linear approximation to the traffic assignment and use that as exogenous traffic assignment information $G$ in the flow estimation model in Program \ref{eq:flow_opt}. One such example is given in Figure \ref{fig:G_example}. In this example, we show 10 OD pairs, and they correspond to 40 paths combined. The values of the matrix represent the proportion of the trips taking the specific path connecting the OD pair. For example, all trips for the OD pair 2 use path 18. For OD pair 1, $20\%$ of trips use paths 0, 15 and 36, respectively, and $40\%$ of trips use path 13.


\begin{table}[!htbp]
    \caption{Time of Day Intervals}\label{tab:tod_intervals}
    \begin{center}
        \begin{tabular}{l l l l l}
        \hline
            TOD Interval & TOD Definition & Time start & Time end & Length (Hours) \\\hline
            1 & AM early & 0:00 & 7:00& 7 \\
            2 & AM peak & 7:00 & 10:00& 3 \\
            3 & Midday & 10:00 & 15:00 & 5 \\
            4 & PM peak & 15:00 & 19:00 & 4 \\
            5 & PM late & 19:00 & 22:00& 3 \\
            6 & Night & 22:00 & 24:00 & 2 \\\hline
        \end{tabular}
    \end{center}
\end{table}

\subsection{Calibration results}
\begin{figure}[!htbp]
  \centering
  \includegraphics[width=0.6\textwidth]{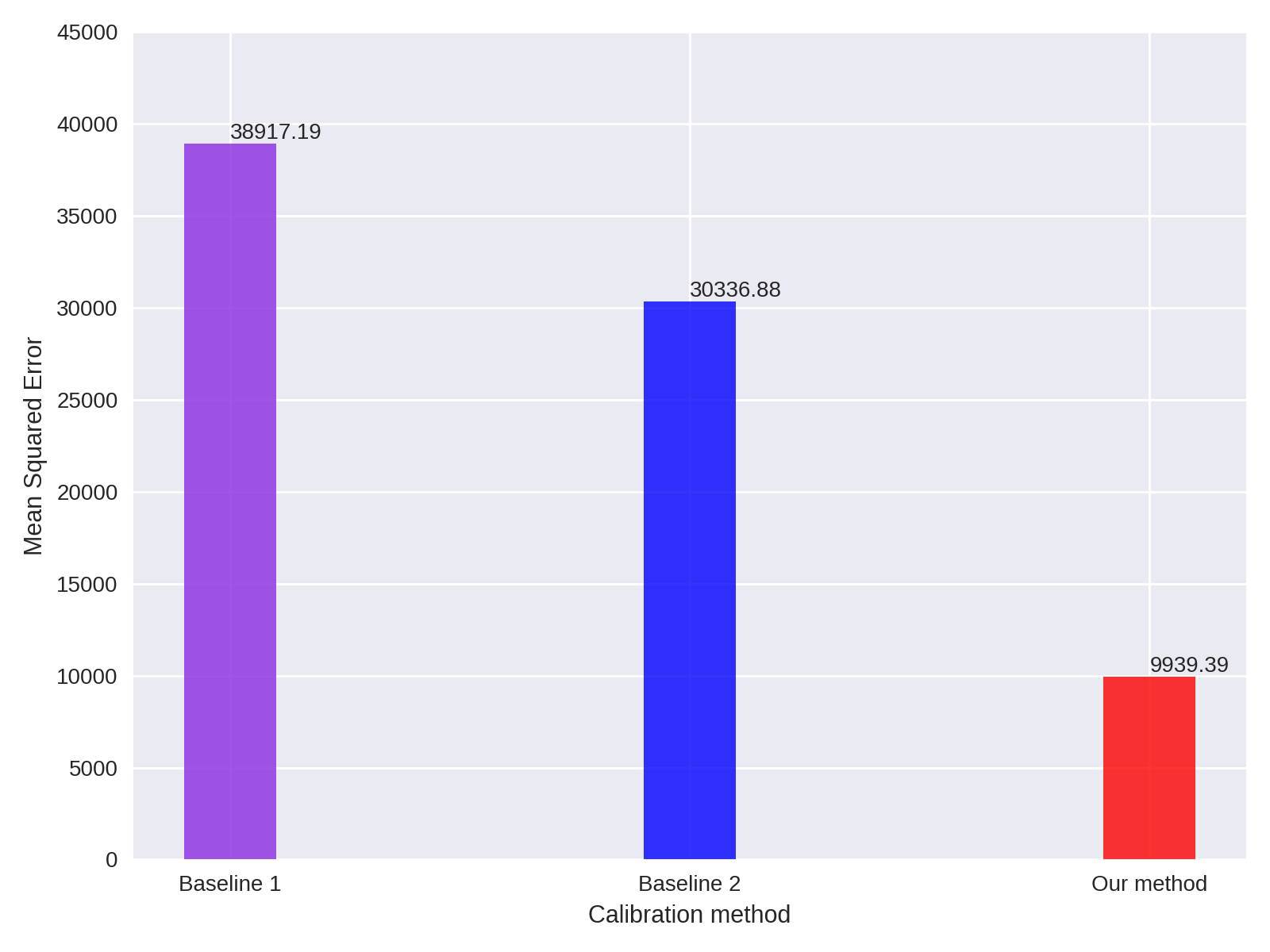}
  \caption{MSE }\label{fig:mse_comparison}
\end{figure}

In this section, we present the aggregated statistics and experiment setup. For comparison purposes, in addition to our proposed calibration method, we also implemented two baseline test cases for the Birmingham network simulation. For baseline 1, we uniformly up-sample the trajectory observations, and adjust the capacity parameters to accommodate as many trips as possible in the network. For baseline 2, we uniformly up-sample the trajectory observations and use the same calibrated capacity parameters obtained from SPSA. The comparison is summarized in Table \ref{table:comparison}. The MSE results in terms of trip travel time (in seconds) are included in Figure \ref{fig:mse_comparison}. Note that for baseline 1, we can only load $90\%$ trips onto the network, even with adjusted parameters for higher link and node capacity. For baseline 2, with the up-sampled trajectory and the same link and node capacity parameter as our calibrated parameters, we can only load $72\%$ trips. Our method demonstrates the improved performance in MSE, and we could load $100\%$ full-scale trips.

\begin{table}[!htbp]
\caption{Summary of performance comparison}
\resizebox{\textwidth}{!}{%
\begin{tabular}{lllll}
\hline
Scenario   & Network flow                  & Link and node capacity                   & System throughput (trips) & MSE (travel time) \\
\hline
Baseline 1 & Up-sampling         & Adjusted parameters with higher capacity & $90\%$                        & 38917.19        \\
Baseline 2 & Up-sampling         & Calibrated parameters                    & $72\%$                          & 30336.88\\
Our method & Network flow estimation & Calibrated parameters                    & $100\%$                         & 9939.39       \\
\hline
\end{tabular}
}
\label{table:comparison}
\end{table}

\begin{figure}[!htbp]
  \centering
  \begin{subfigure}[b]{0.31\textwidth}
    \centering
    \includegraphics[width=\textwidth]{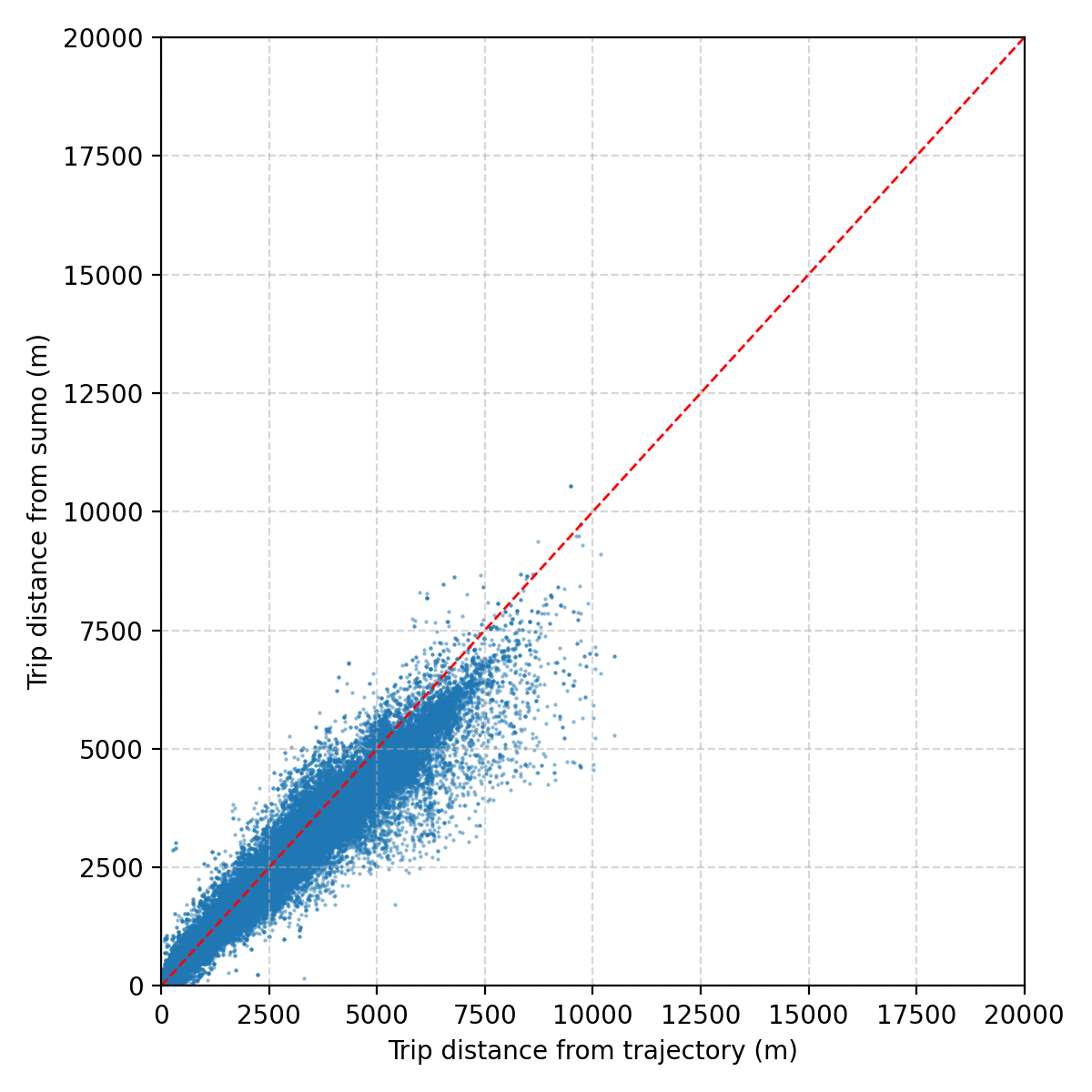}
    \caption{Trip distance}\label{fig:baseline1_1}
  \end{subfigure}
  \hfill
  \begin{subfigure}[b]{0.31\textwidth}
    \centering
    \includegraphics[width=\textwidth]{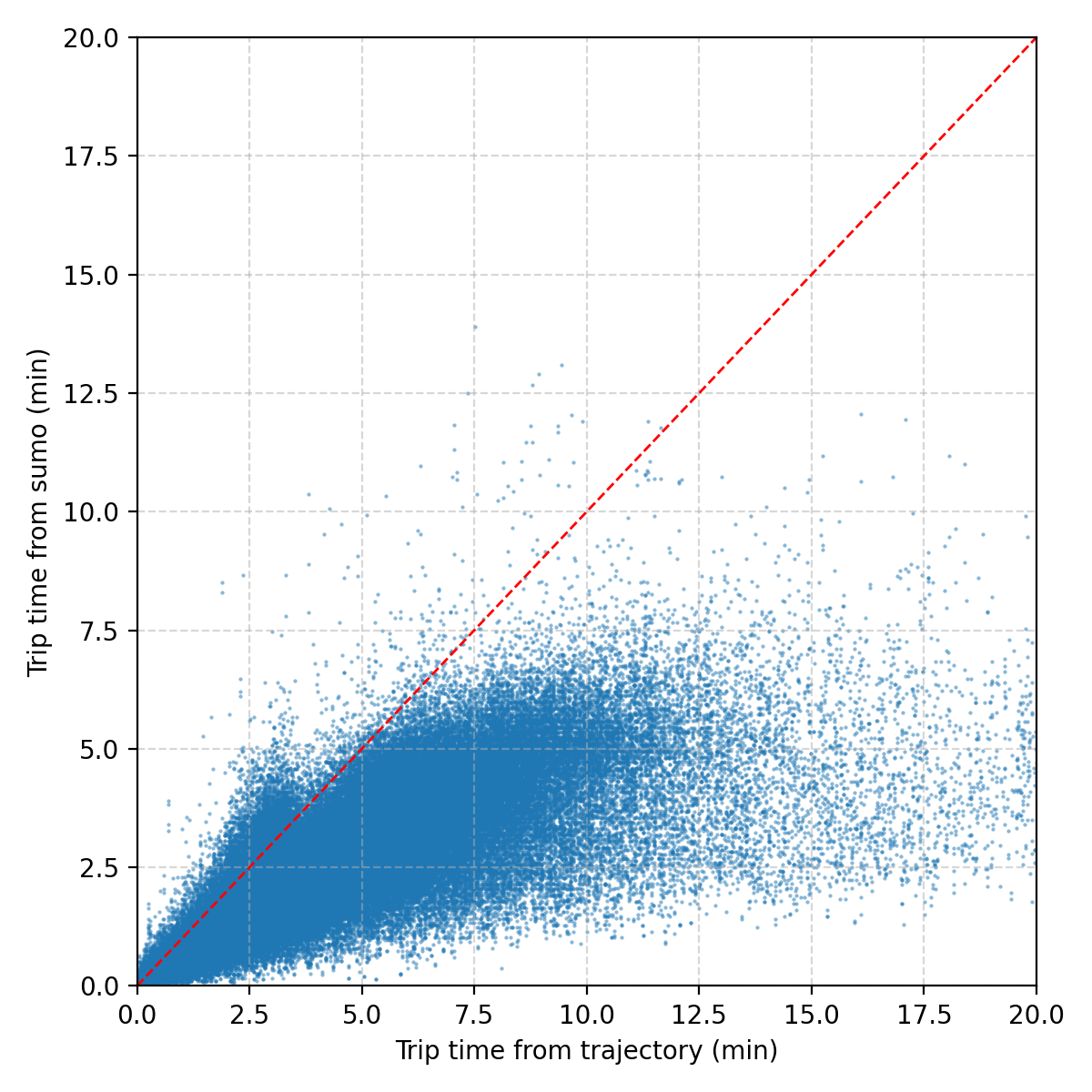}
    \caption{Trip travel time}\label{fig:baseline1_2}
  \end{subfigure}
  \hfill
  \begin{subfigure}[b]{0.31\textwidth}
    \centering
    \includegraphics[width=\textwidth]{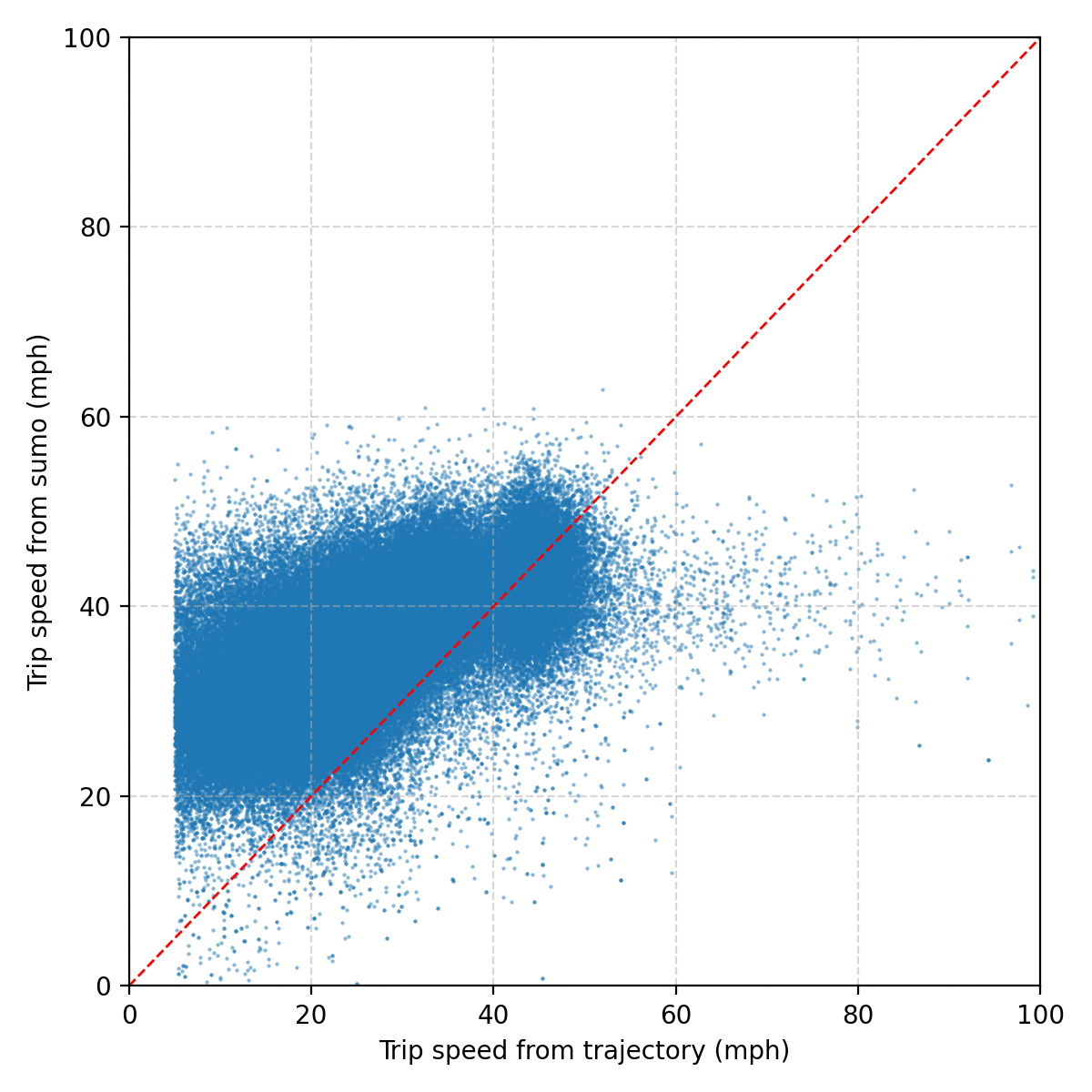}
    \caption{Trip travel speed}\label{fig:baseline1_3}
  \end{subfigure}
  \caption{ Performance of baseline 1 with maximum $90\%$ of total trips}\label{fig:baseline1}
\end{figure}

\begin{figure}[!htbp]
  \centering
  \begin{subfigure}[b]{0.31\textwidth}
    \centering
    \includegraphics[width=\textwidth]{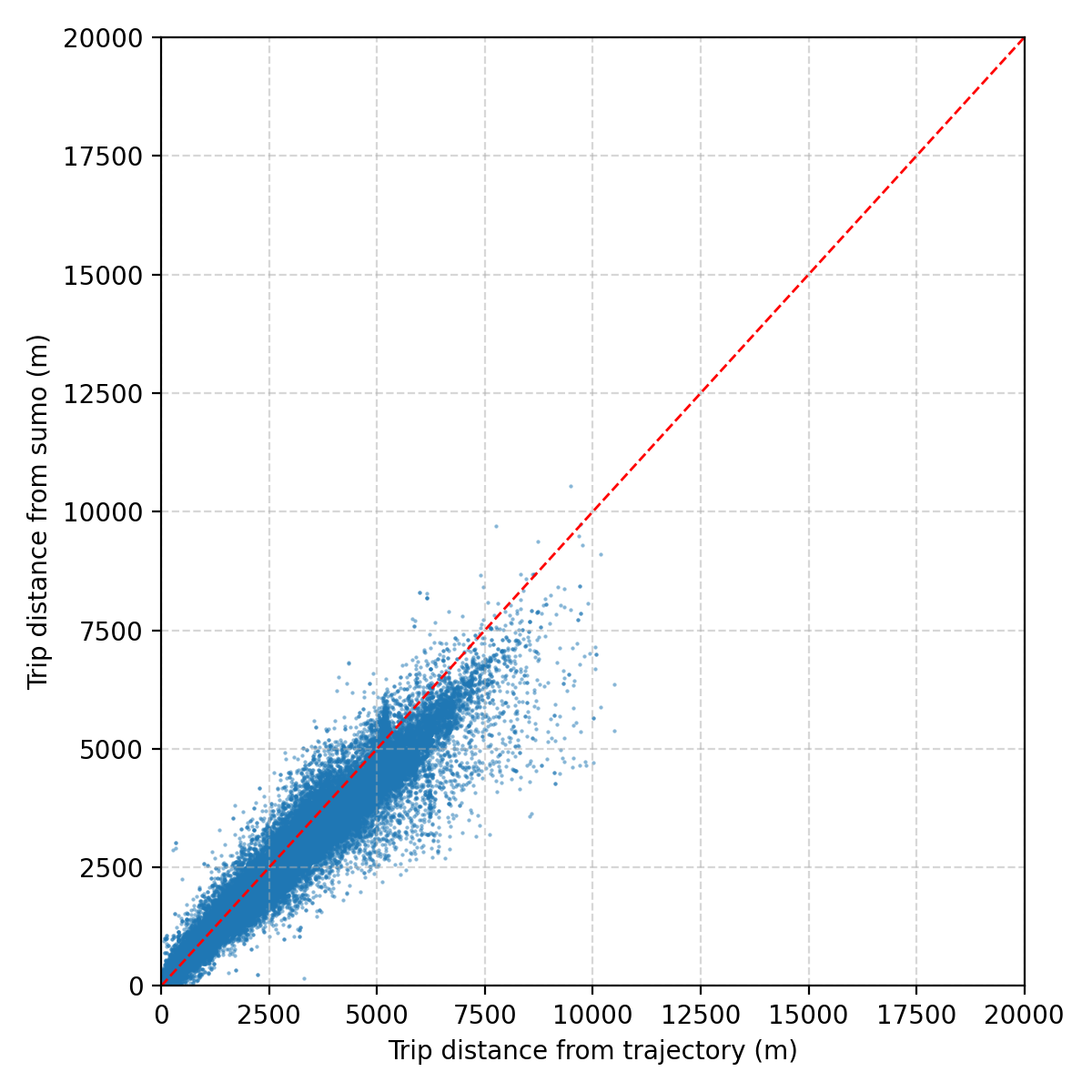}
    \caption{Trip distance}\label{fig:baseline2_1}
  \end{subfigure}
  \hfill
  \begin{subfigure}[b]{0.31\textwidth}
    \centering
    \includegraphics[width=\textwidth]{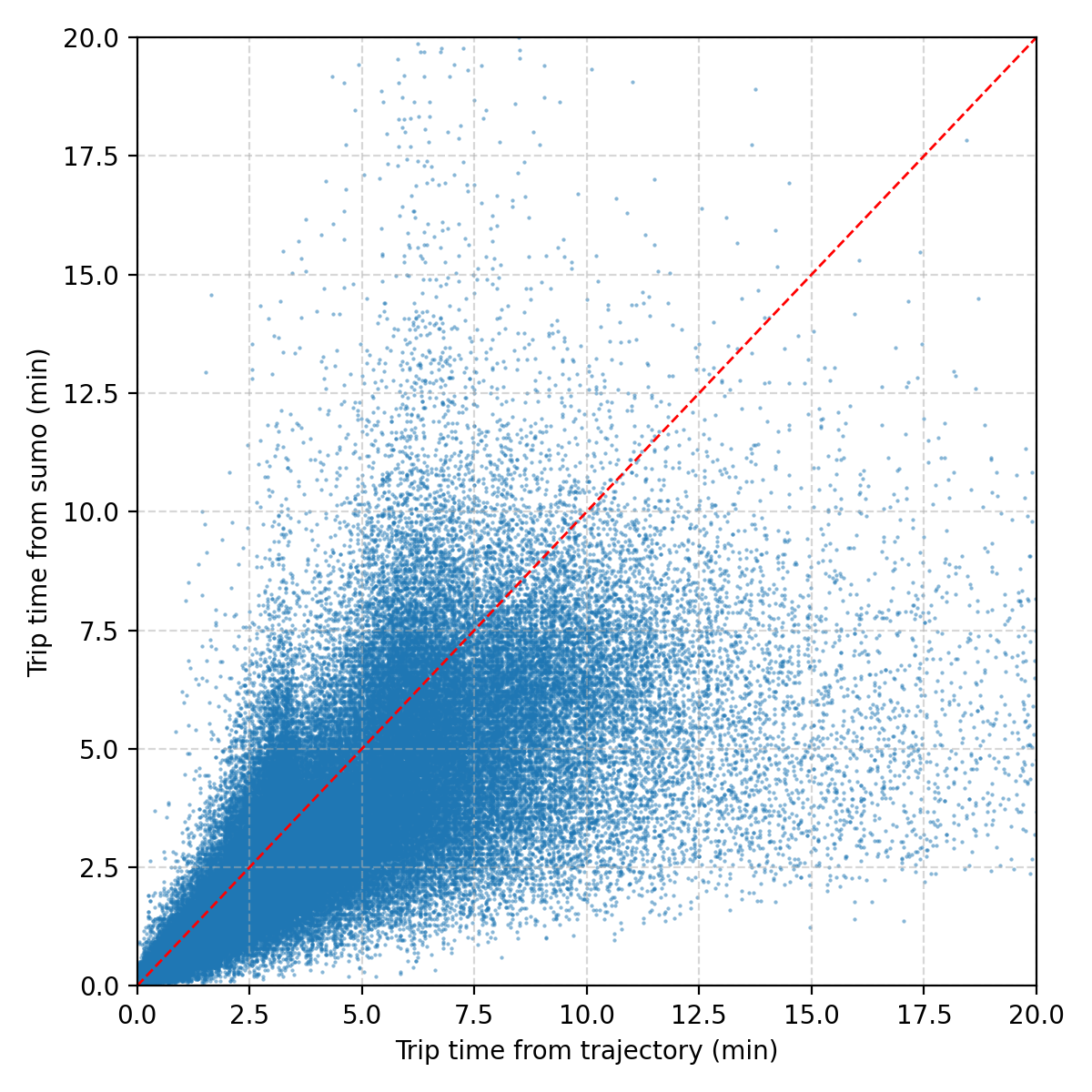}
    \caption{Trip travel time}\label{fig:baseline2_2}
  \end{subfigure}
  \hfill
  \begin{subfigure}[b]{0.31\textwidth}
    \centering
    \includegraphics[width=\textwidth]{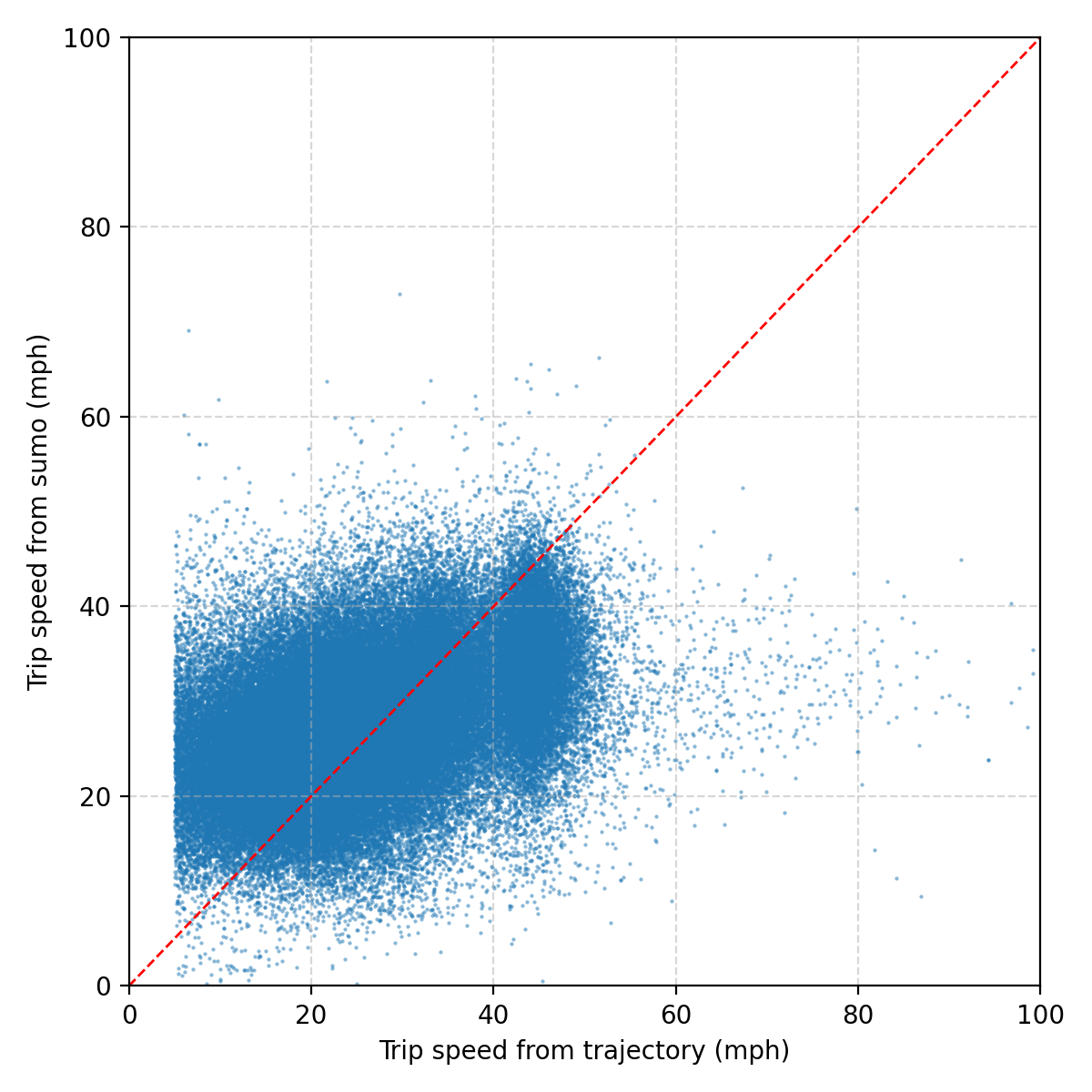}
    \caption{Trip travel speed}\label{fig:baseline2_3}
  \end{subfigure}
  \caption{Performance of baseline 2 with maximum $72\%$ of total trips}\label{fig:baseline2}
\end{figure}

\begin{figure}[!htbp]
  \centering
  \begin{subfigure}[b]{0.31\textwidth}
    \centering
    \includegraphics[width=\textwidth]{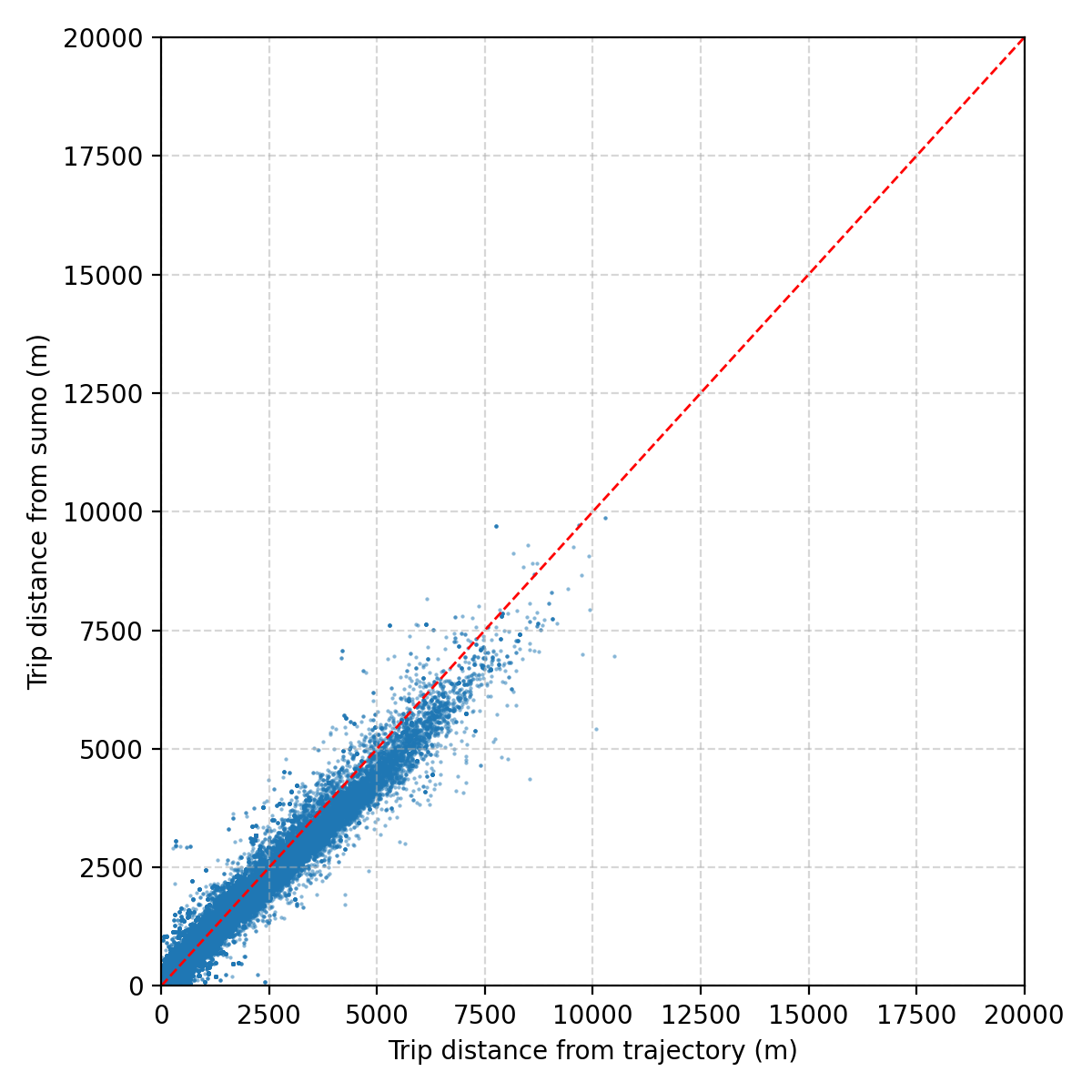}
    \caption{Trip distance}\label{fig:subfig1}
  \end{subfigure}
  \hfill
  \begin{subfigure}[b]{0.31\textwidth}
    \centering
    \includegraphics[width=\textwidth]{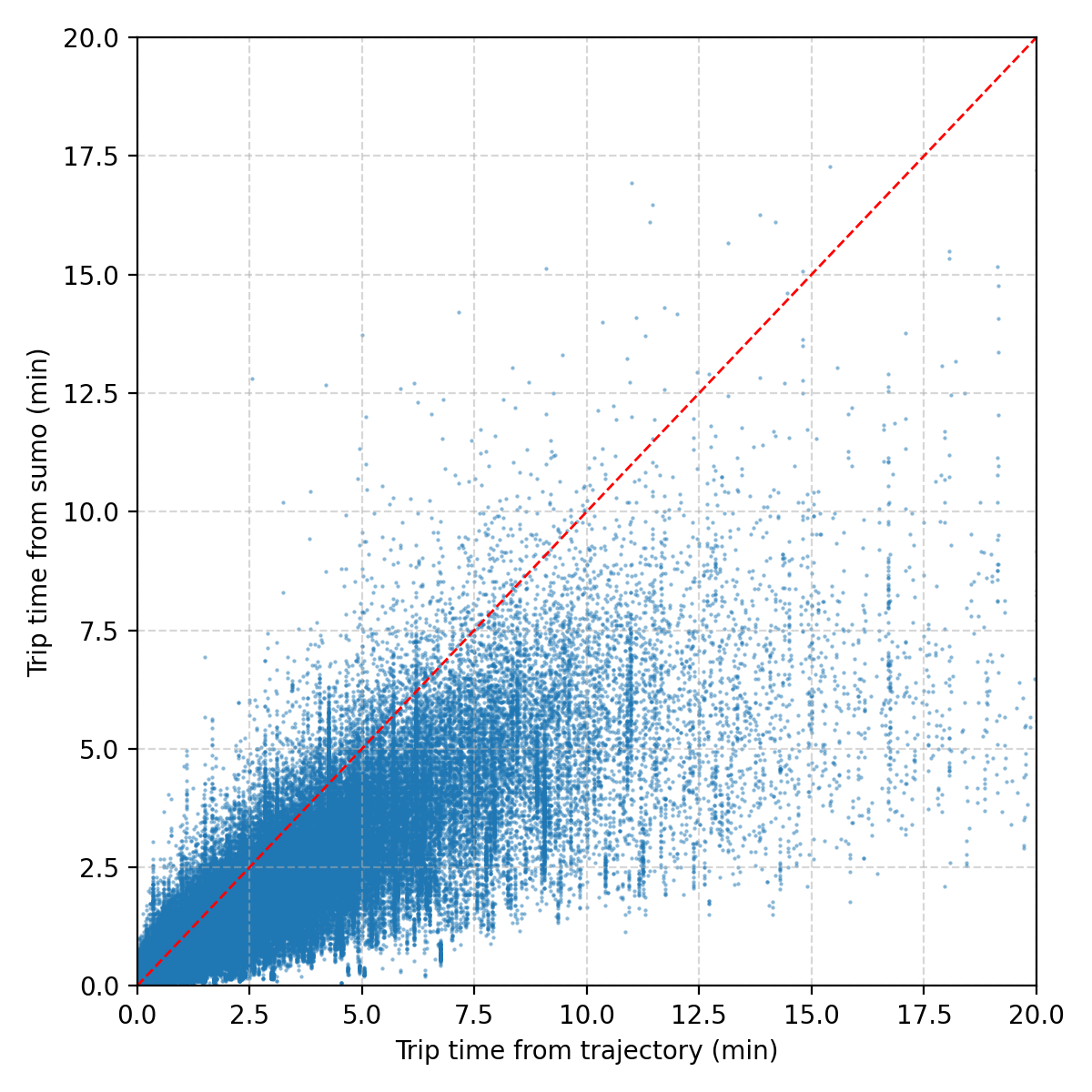}
    \caption{Trip travel time}\label{fig:subfig2}
  \end{subfigure}
  \hfill
  \begin{subfigure}[b]{0.31\textwidth}
    \centering
    \includegraphics[width=\textwidth]{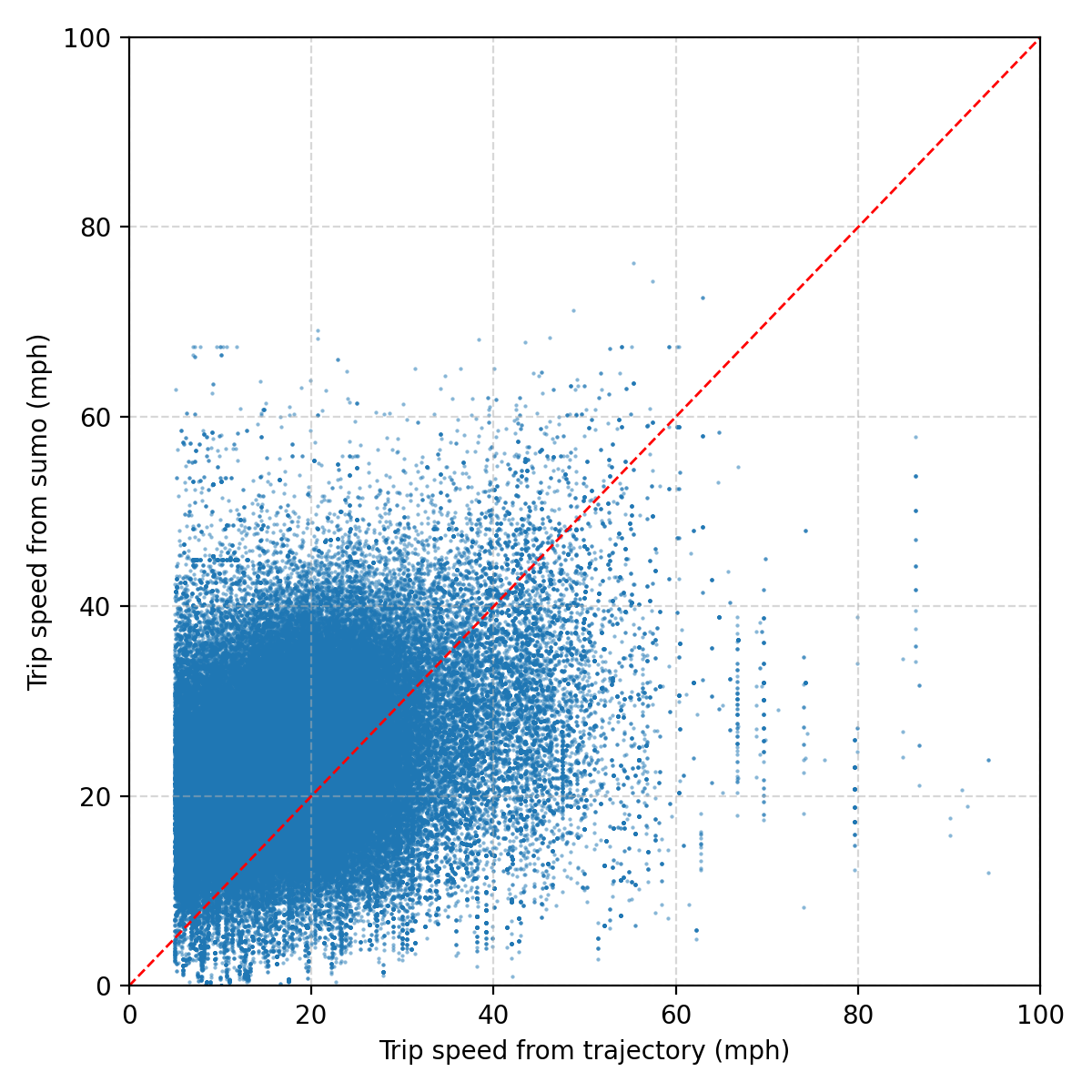}
    \caption{Trip travel speed}\label{fig:subfig3}
  \end{subfigure}
  \caption{Performance of our calibration method with full-scale $100\%$ total trips}\label{fig:our}
\end{figure}

\begin{figure}[!htbp]
  \centering
  \includegraphics[width=0.55\textwidth]{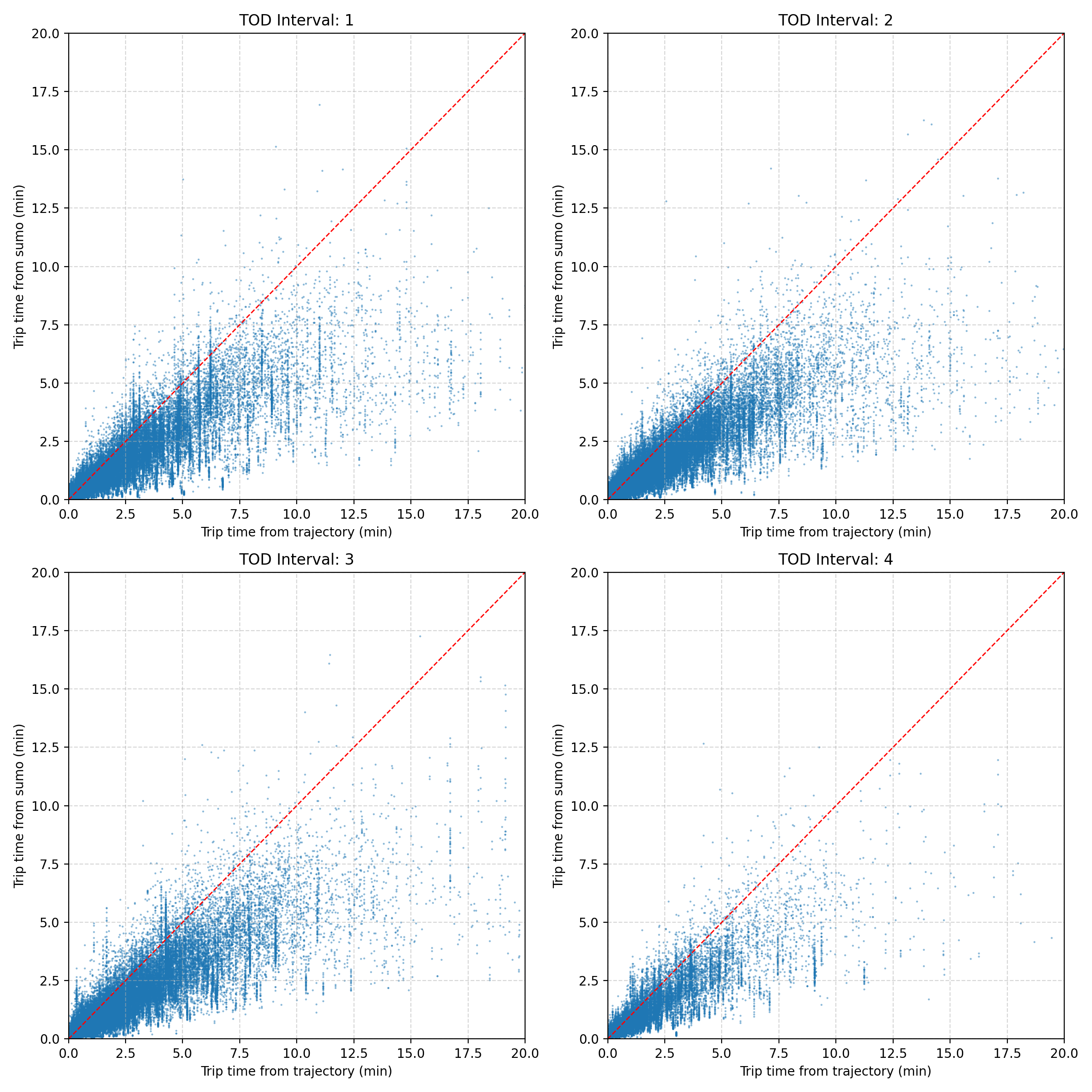}
  \caption{Performance of our calibration method in trip travel time breakdown by TOD}\label{fig:tod_breakdown}
\end{figure}

Next, we present the comparison of trip distance, trip travel time, and trip speed for these three methods. The performance of baseline 1, baseline 2, and our method is shown in Figure \ref{fig:baseline1}, Figure \ref{fig:baseline2} and Figure \ref{fig:our}, respectively. Among the three cases, our method produces the most consistent trip distance compared to trajectory observations. In baseline 1, the simulated trip speeds from SUMO are generally larger than those from the trajectory. However, the network can only support $90\%$ of total trips. 
In baseline 2, the trip speeds are closer to the observed trajectory. However, there exist many extremely fast trips and slow trips, as shown off-diagonal in Figure \ref{fig:baseline2_2}. In addition, the system throughput is reduced to $72\%$. 
In our method, the simulated trip travel times are slightly biased towards smaller values, implying less congestion or faster drivers compared to the real-world cases. However, the MSE of $9939.39$ from our method is much better compared to $38917.19$ from baseline 1 and $30336.88$ from baseline 2. In addition, the trip travel time performance is broken down by 4 main TOD intervals, as shown in Figure \ref{fig:tod_breakdown}. The performance is consistent across different congestion levels during peak and off-peak time intervals. 

One main reason that the baseline cases can support fewer trips and result in larger travel time errors is due to the unreasonable route choice. The long trips can be over-sampled compared to the ground truth during the uniform up-sampling process for the baseline cases. In reality, the longer paths connecting an OD pair may only be used during abnormal traffic congestion patterns or due to personal habitual preferences. However, the uniform up-sampling cannot differentiate the rational route choices versus the abnormal conditions. In our calibration framework, this over-penalization is mitigated by the presence of a network flow estimation model. This makes our method more robust and realistic compared to the real-world situations.

\section{Discussion and Conclusion}
This paper presents a fast automatic calibration framework based on vehicle trajectory data for mesoscopic traffic simulation models. The proposed method provides a systematic approach to jointly calibrating demand and supply components for a network-level simulation. The proposed method is demonstrated through a case study of the City of Birmingham in Michigan. Our results demonstrate superior performance compared with two baseline models in trip-level performance metrics. 

A future opportunity is the development of an iterative procedure that simultaneously combines the demand stage and the supply stage. Currently, we have established a network flow estimation model to approximate the network equilibrium state based on trajectory data for multiple TOD intervals. Then the capacity and driving behavior parameters are calibrated thereafter for the same approximated network equilibrium. A worthy step is to combine the demand step and the supply step shown in this paper and design an iterative procedure to better approximate the global optimal for the system state. Another future direction is the acceleration of the SUMO simulation regarding incremental updates in the demand and supply input throughout the iterations. Related to the design of online and incremental optimization, the calibration process would greatly benefit from the computational complexity savings from the accelerated processes. By exploring these directions, we can further improve the efficiency and accuracy of traffic simulation models, ultimately leading to better-informed transportation planning and management.

\bibliographystyle{apalike}  
\bibliography{references}  

\end{document}